\documentclass{gtart}

\def\ifplaintex{\expandafter\ifx\csname documentclass\endcsname\relax}

\ifplaintex 
\else
\fi

\def\gtm{{\mathsurround=0pt\it $\cal G\mskip-2mu$eometry \&\ 
$\cal T\!\!$opology $\cal M\mskip-1mu$onographs}}    

\def\gtp{{\mathsurround=0pt\it $\cal G\mskip-2mu$eometry \&\ 
$\cal T\!\!$opology $\cal P\!$ublications}}  

\def\recd{{\small Received:\qua\receiveddate\ifx\reviseddate\relax
\else\qquad Revised:\qua\reviseddate\fi\par}} 

\def\volumenumber#1{\def\thevolumenumber{#1}}
\def\volumeyear#1{\def\thevolumeyear{#1}}
\def\volumename#1{\def\thevolumename{#1}}
\def\papernumber#1{\def\thepapernumber{#1}}
\def\pagenumbers#1#2{\def\startpage{#1}\def\finishpage{#2}}
\def\published#1{\def\publishdate{#1}}
\def\received#1{\def\receiveddate{#1}}
\def\revised#1{\def\reviseddate{#1}}
\def\accepted#1{\def\accepteddate{#1}}

\long\def\asciiabstract#1{\long\def\theasciiabstract{#1}}

\let\\\par
\let\thevolumenumber\relax\let\thepapernumber\relax
\let\thevolumeyear\relax\let\startpage\relax
\let\finishpage\relax\let\publishdate\relax\let\receiveddate\relax
\let\reviseddate\relax\let\accepteddate\relax\let\theasciititle\relax
\let\theasciiauthors\relax
\let\theasciiabstract\relax

\let\theerratum\relax\let\theasciiemail\relax
\let\theshortauthors\relax\let\theshorttitle\relax

\def\startpage{1}\def\finishpage{15}\def\thepapernumber{77}

\volumenumber{2}
\volumename{Proceedings of the Kirbyfest}
\volumeyear{1999}

\long\def\maketitlep{   

\count0=\startpage

\gtm\nl        
{\small Volume \thevolumenumber: \thevolumename\nl 
\ifx\theerratum\relax\else Erratum \erratumnumber\nl\fi
Pages \startpage--\finishpage\nl}

\vglue 0.1truein   

{\parskip=0pt\leftskip 0pt plus 1fil\def\\{\par\smallskip}{\ifplaintex\large
\else\Large\fi\bf\thetitle}\par\medskip}   
\vglue 0.05truein 

{\parskip=0pt\leftskip 0pt plus 1fil\def\\{\par}{\sc\theauthors}
\par\medskip}%
 
\vglue 0.03truein 

{\small\leftskip 25pt\rightskip 25pt{\bf Abstract}\stdspace\theabstract

{\bf AMS Classification}\stdspace\theprimaryclass
\ifx\thesecondaryclass\relax\else; \thesecondaryclass\fi\par
{\bf Keywords}\stdspace \thekeywords\par}\vglue 7pt

}   

\font\phead=cmsl9 scaled 950
\font=cmsl9 scaled 1050
\font\pnum=cmbx10 scaled 913
\font\lnum=cmbx10 
\font\pfoot=cmsl9 scaled 950
\font=cmsl9 scaled 1050
\ifplaintex
\headline{\vbox to 0pt{\vskip -4.5mm\line{\small\phead\ifnum
\count0=\startpage ISSN 1464-8997 (on line)
1464-8989 (printed) \hfill {\pnum\folio}\else\ifodd\count0\def\\{ }%
\ifx\theshorttitle\relax\thetitle\else\theshorttitle\fi\hfill{\pnum\folio}
\else\def\\{ and }{\pnum\folio}\hfill\ifx\theshortauthors\relax\theauthors
\else\theshortauthors\fi\fi\fi}\vss}}
\footline{\vbox to 0pt{\vglue 0mm\line{\small\pfoot\ifnum\count0=\startpage
Published \publishdate:\qua\copyright\ \gtp\hfill\else
\gtm, Volume \thevolumenumber\ (\thevolumeyear)\hfill\fi}\vss
}}
\else
\makeatletter
\def\@oddhead{{\small\lhead\ifnum\count0=\startpage ISSN 1464-8997 (on line)
1464-8989 (printed) \hfill {\lnum\number\count0}\else\ifodd\count0
\def\\{ }\ifx\theshorttitle\relax \thetitle \else\theshorttitle\fi\hfill
{\lnum\number\count0}\else\def\\{ and }{\lnum\number\count0}
\hfill\ifx\theshortauthors\relax 
\theauthors\else\theshortauthors\fi\fi\fi}}\def\@evenhead{@oddhead}
\def\@oddfoot{\small\lfoot\ifnum\count0=\startpage Published \publishdate:\qua\copyright\ \gtp\hfill\else
\gtm, Volume \thevolumenumber\ (\thevolumeyear)\hfill\fi}
\def\@evenfoot{@oddfoot}
\makeatother
\fi

\let\maketitlepage\maketitlep
\let\makeshorttitle\maketitlepage
\let\maketitle\maketitlepage

\newwrite\gtoutfile
\long\gdef\makeheadfile{  
{\def\\{, }\def\s{ }
\immediate\openout\gtoutfile head.xxx
\immediate\write\gtoutfile{To: math@arxiv.org}
\immediate\write\gtoutfile{Subject: put OR rep NNNNN:ppppp}
\immediate\write\gtoutfile{--text follows this line--}
\immediate\write\gtoutfile{Proxy-for: \ifx\theasciiauthors\relax
\theauthors\else\theasciiauthors\fi\s<\ifx\theasciiemail\relax\theemail\else\theasciiemail\fi>}
\immediate\write\gtoutfile{\noexpand\\}
\immediate\write\gtoutfile{Authors: \ifx\theasciiauthors\relax
\theauthors\else\theasciiauthors\fi}
{\def\\{ }\immediate\write\gtoutfile{Title: \ifx\theasciititle\relax
\thetitle\else\theasciititle\fi}}
\immediate\write\gtoutfile{Subj-class: GT or SG, GR etc}
\immediate\write\gtoutfile{MSC-class: \theprimaryclass\ifx\thesecondaryclass\relax\else, \thesecondaryclass\fi}
\immediate\write\gtoutfile{Journal-ref: Geom. Topol. Monogr. \thevolumenumber\s
(\thevolumeyear) \startpage-\finishpage}
\immediate\write\gtoutfile{Comments: Published by Geometry and Topology Monographs at}
\immediate\write\gtoutfile{\s\s\s  http://www.maths.warwick.ac.uk/gt/GTMon\thevolumenumber/paper\thepapernumber.abs.html}
\immediate\write\gtoutfile{\noexpand\\}
\immediate\write\gtoutfile{}
\ifx\theasciiabstract\relax
\immediate\write\gtoutfile{\theabstract}\else
\immediate\write\gtoutfile{\theasciiabstract}\fi
\immediate\write\gtoutfile{}
\immediate\write\gtoutfile{\noexpand\\}
\immediate\write\gtoutfile{}
\immediate\closeout\gtoutfile}}  

\def\maketitlepage{\maketitlep\makeheadfile}
\let\makeshorttitle\maketitlepage
\let\maketitle\maketitlepage

\volumenumber{4}
\volumename{Invariants of knots and 3-manifolds (Kyoto 2001)}
\volumeyear{2002}
\papernumber{10}
\pagenumbers{143}{160}
\received{19 December 2001}
\revised{9 April 2002}
\accepted{22 July 2002}
\published{19 September 2002}

\usepackage{amssymb,amsmath,graphicx}

\newtheorem{thm}{Theorem}[section]   
\newtheorem{lem}[thm]{Lemma}         
\newtheorem{prop}[thm]{Proposition}

\theoremstyle{definition}

\begin{document}
\title[Loop spaces of configuration spaces]
{Loop spaces of configuration spaces\\and finite type invariants}

\authors{Toshitake Kohno}
\address{Graduate School of Mathematical Sciences\\University of
Tokyo, Tokyo 153-8914 Japan}
\email{kohno@ms.u-tokyo.ac.jp}

\begin{abstract}   
The total homology of the loop space of the configuration space
of ordered distinct $n$ points in ${\bf R}^m$ has a structure of
a Hopf algebra defined by the 4-term relations if $m \geq 3$.
We describe a relation of between the cohomology of this
loop space and the set of finite type invariants for
the pure braid group with $n$ strands.
Based on this we give
expressions of certain link invariants as integrals over
cycles of the above loop space.
\end{abstract}

\asciiabstract{
The total homology of the loop space of the configuration space of
ordered distinct n points in R^m has a structure of a Hopf algebra
defined by the 4-term relations if m>2.  We describe a relation of
between the cohomology of this loop space and the set of finite type
invariants for the pure braid group with n strands.  Based on this we
give expressions of certain link invariants as integrals over cycles
of the above loop space.}

\maketitle

\cl{\small\em Dedicated to Professor Mitsuyoshi Kato on his
sixtieth birthday}

\section{Introduction}

The purpose of this article is to present a new approach
for finite type invariants of braids based on
the loop spaces of configuration spaces.
For a smooth manifold $M$ with a base point $x_0$ we denote by
$\Omega M$ the loop space consisting of the
piecewise smooth loops
$\gamma : I \rightarrow M$ with
$\gamma (0)= \gamma (1) =x_0$.
The main object is the loop space
$\Omega {\rm Conf}_n({\bf R}^m)$
where
${\rm Conf}_n({\bf R}^m)$ stands for the configuration
space of ordered distinct $n$ points in ${\bf R}^m$.
It follows from the work of Cohen and Gitler \cite{CG}
that, in the case $m\geq 3$,
the total homology of the loop space
$\Omega {\rm Conf}_n({\bf R}^m)$
has a structure of the
universal enveloping algebra of a Lie algebra
defined by 4-term relations in the graded sense.
Based on this one can show that the space of
finite type invariants for pure braids
is naturally isomorphic to the total cohomology
of ${\rm Conf}_n({\bf R}^m)$
if $m \geq 3$.
The total homology of the loop space
$\Omega {\rm Conf}_n({\bf R}^m)$
has a uniform structure for
$m \geq 3$ except a shift of grading.
In this article, we mainly deal with the
case $m=3$.
 From our point of view the space of weight systems
for Vassiliev invariants
for knots in the sense of \cite{V}
lies in the direct sum
$\bigoplus_{n \geq 2} H^*(\Omega {\rm Conf}_n({\bf R}^m))$.
The de Rham cohomology of
$\Omega {\rm Conf}_n({\bf R}^m)$ is described by
Chen's iterated integrals of Green forms.
We show that a certain link invariant can be
expressed as an integral of a differential form
on $\Omega {\rm Conf}_n({\bf R}^m)$
over a cycle defined associated with a link.
It is still work in progress to obtain a
general formula for any finite type invariant.
Let us recall that
based on the work of Guadagnini, Martellini and
Mintchev \cite{GMM} for Chern-Simons perturbative
theory, Bott and Taubes developed a method
to express a knot
invariant by means of iterated integrals of
Green forms associated with Feynman diagrams.
In this approach one needs to add integrals
associated with graphs with trivalent
vertices to obtain a knot invariant.
Our method does not use trivalent vertices, but
we add
integrals of some non closed
differential forms on the loop
spaces of
configuration spaces to obtain topological invariants.

The paper is organized in the following way.
In Section 2, we give a description of the
homology of the loop space
$\Omega {\rm Conf}_n({\bf R}^m)$
by means of the graded 4-term relations.
In Section 3, we briefly review Chen's
iterated integrals describing the de Rham
cohomology of loop spaces.
In Section 4,
we discuss a relation between
the above homology and finite type invariants
for braids.
In Section 5, we express certain link invariants
in terms of an integral of a differential form
on the loop space of the configuration spaces.
Section 6 is devoted to a brief overview
of a work in progress in \cite{CKX}
on the homology of the loop space of
orbit configuration spaces associated to
Fuchsian groups.
A more detailed account of this subject
will appear
elsewhere.

\section
{Homology of the loop spaces}

We denote by ${\rm Conf}_n(X)$
the configuration space of ordered distinct $n$ points
in a space $X$. Namely, we set
$$
{\rm Conf}_n(X) =
\{ (x_1, \cdots, x_n) \in X^n  \ ;  \
x_i \neq x_j  \ \mbox{if} \ i \neq j
\}.
$$
Let $\Delta_{ij}$ be the diagonal set
of $X^n$ defined by $x_i = x_j$.
We will deal with the loop space
$\Omega {\rm Conf}_n(X)$
of the configuration space ${\rm Conf}_n(X)$
with a fixed base
point $x_0$.

For the rest of this section we consider the case
$X={\bf R}^m$ with $m\geq 3$.
The boundary of a tubular neighbourhood
of the diagonal set $\Delta_{ij}$
is identified with the tangent
sphere bundle of $X$ and we have a map
$$
\gamma_{ij} :  S^{m-1} \rightarrow {\rm Conf}_n({\bf R}^m),  \
1 \leq i < j \leq n,
$$
so that $\gamma_{ij}(S^{m-1})$ has the linking number
$1$ with the diagonal set $\Delta_{ij}$.
Since
$S^{m-1}$ is considered to be
the suspension of its
equator
$S^{m-2}$, we have a natural map
$$
\alpha_{ij} :  S^{m-2} \rightarrow \Omega {\rm Conf}_n({\bf R}^m),  \
1 \leq i < j \leq n,
$$
induced from $\gamma_{ij}$.

First we describe the cohomology ring of the configuration
space ${\rm Conf}_n({\bf R}^m)$.
Let $\omega$ be a homogeneous $(m-1)$-form on
${\bf R}^m \setminus \{ {\bf 0}\}$ defining a standard volume
form of the unit sphere
$S^{m-1}$.
Namely, $\omega$ satisfies
\begin{align*}
&\int_{S^{m-1}}\omega =1,  \\
&\omega (\lambda {\bf x})=({\rm sgn} \ \lambda)^m \omega ({\bf x})
\end{align*}
for ${\bf x} \in
{\bf R}^m \setminus \{ {\bf 0}\}$ and
$\lambda \in {\bf R} \setminus \{ {\bf 0}\}$.
Here
${\rm sgn} \ \lambda$ stands for $\lambda / |\lambda|$.
We define the Green form
$\omega_{ij}$ by
$$
\omega_{ij}({\bf x}_1, \cdots, {\bf x}_n)
= \omega ({\bf x}_j - {\bf x}_i),
\quad {\bf x}_1, \cdots, {\bf x}_n \in {\bf R}^m,
$$
which is an $(m-1)$-form on ${\rm Conf}_n({\bf R}^m)$
for $i \neq j$.
The de Rham cohomology class $[\omega_{ij}]$
defines an integral cohomology
class and is dual to the homology class
$[\gamma_{ij}]$.
We put ${\xi}_{ij}=[\omega_{ij}]$. It is known that the
cohomology ring $H^*({\rm Conf}_n({\bf R}^m) ; {\bf Z})$ is
  generated by
${\xi}_{ij}$, $1 \leq i \neq j \leq n$,
  with relations
\begin{align*}
& {\xi}_{ij}^2 = 0 \\
& {\xi}_{ij} = (-1)^m {\xi}_{ji} \\
& {\xi}_{ij}{\xi}_{jk} + {\xi}_{jk}{\xi}_{ik} + {\xi}_{ik}{\xi}_{ij}=0, \ i<j<k
\end{align*}
where ${\rm deg} \ \omega_{ij}=m-1$ (see \cite{Cohen}).

In general,
the total homology
$$
H_*(\Omega M ; {\bf Z})
=
\bigoplus_{j\geq 0}
H_j(\Omega M ; {\bf Z})
$$
of the loop space of $M$
is equipped with a product
$$
H_i(\Omega M ; {\bf Z})
\otimes
H_j(\Omega M ; {\bf Z})
\rightarrow
H_{i+j}(\Omega M ; {\bf Z})
$$
induced from the composition of loops.
We also have a coproduct
$$
H_k(\Omega M ; {\bf Z})
\rightarrow
\oplus_{i+j=k}
H_i(\Omega M ; {\bf Z})
\otimes
H_j(\Omega M ; {\bf Z})
$$
induced from the cup product homomorphism on
cochains.
Let us now investigate the total homology of the
loop space
$\Omega {\rm Conf}_n({\bf R}^m)$ in the case
$m \geq 3$
as a Hopf algebra with the above product and
coproduct.

Let us first consider the simplest example
$\Omega {\rm Conf}_2({\bf R}^m)$.
Here the configuration space
${\rm Conf}_2({\bf R}^m)$ is homotopy equivalent to
$S^{m-1}$. The structure of the total homology of the
loop space of a sphere was determined by
Bott and Samelson \cite{BS} (see also \cite{Chen}).
We have isomorphisms of Hopf algebras
$$
H_*(\Omega {\rm Conf}_2({\bf R}^m) ; {\bf Z}) \cong
H_*(\Omega S^{m-2} ; {\bf Z}) \cong
{\bf Z}[X_{12}]
$$
where ${\bf Z}[X_{12}]$ stands for the polynomial algebra
with one indeterminate with
$\deg X_{12}=m-2$ and $X_{12}$ corresponds to
the
homology class represented by $\alpha_{12}$ defined
above.


In general we have relations among the homology classes
$X_{ij}=[\alpha_{ij}]$ analogous to the 4-term relations.
We set $X_{ij}=(-1)^{m-2}X_{ji}$ for $i>j$.

\begin{prop}
In $H_*(\Omega {\rm Conf}_n({\bf R}^m) ; {\bf Z})$,
the homology classes $X_{ij}$ satisfy the relation
$$
[X_{ij},  X_{ik} + X_{jk}]=0
$$
for $i<j<k$. Here $\deg X_{ij}=m-2$ and we consider
the Lie bracket in the graded sense.
\end{prop}

\begin{proof}
We fix $i, j$ and $k$ with $1 \leq i < j < k \leq n$
and define
$$
\varphi : S^{m-1} \times S^{m-1}
\rightarrow
{\rm Conf}_n({\bf R}^m)
$$
in the following way.
We take ${\bf u} \in {\bf R}^m$ with $\| {\bf u} \|=1$ and fix
$r_1$ and $r_2$ such that $0<r_1<r_2<1$.
For ${\bf x}_1, {\bf x}_2 \in S^{m-2}$ and $1\leq l \leq n$
we set
\begin{equation*}
\varphi_{l}({\bf x}_1, {\bf x}_2) =
\begin{cases}
{l}{\bf u},  & \quad l \neq j, k \\
i{\bf u} + r_1 {\bf x}_1, & \quad {l}=j \\
i{\bf u} + r_2 {\bf x}_2, & \quad {l}=k
\end{cases}
\end{equation*}
and $\varphi$ is defined to be
$\varphi ({\bf x}_1, {\bf x}_2) = (\varphi_1 ({\bf x}_1, {\bf x}_2),
\cdots, \varphi_{n}({\bf x}_1, {\bf x}_2))$.
We denote by $\alpha \in H_{m-1}(S^{m-1}; {\bf Z})$
the fundamental homology class.
Let us notice that $\alpha$ determines a generator
$\widehat{\alpha}$ of
the total homology of the loop space
$H_*(\Omega S^{m-1}; {\bf Z})$.
We consider the induced homomorphism
$$
\varphi_* :
H_{m-1}(S^{m-1} \times S^{m-1} ; {\bf Z}) \rightarrow
H_{m-1}({\rm Conf}_n({\bf R}^m) ; {\bf Z}).
$$
Then we have
$$
\langle {\xi}_{ij}, \varphi_* (\alpha \times 1) \rangle =1, \
\langle {\xi}_{ik}, \varphi_* (\alpha \times 1) \rangle =0, \
\langle {\xi}_{jk}, \varphi_* (\alpha \times 1) \rangle =0
$$
and
$$
\langle {\xi}_{ij}, \varphi_* (1 \times \alpha) \rangle =0, \
\langle {\xi}_{ik}, \varphi_* (1 \times \alpha) \rangle =1, \
\langle {\xi}_{jk}, \varphi_* (1 \times \alpha) \rangle =1,
$$
which implies
$$
\varphi_* (\alpha \times 1) = \alpha_{ij}, \quad
\varphi_* (1 \times \alpha) = \alpha_{ik} + \alpha_{jk}.
$$
The map $\varphi$ gives
$$
\Omega\varphi :
\Omega (S^{m-1} \times S^{m-1}) \rightarrow
\Omega {\rm Conf}_n({\bf R}^m).
$$
Let us consider the induced homomorphism
$$
\Omega\varphi_* :
H_*(\Omega (S^{m-1} \times S^{m-1}) ; {\bf Z}) \rightarrow
H_*(\Omega {\rm Conf}_n({\bf R}^m) ; {\bf Z}).
$$
The homology classes $\alpha \times 1$ and $1 \times \alpha$
determine the homology classes, say $\widehat{\alpha}_1$ and
$\widehat{\alpha}_2$, in
$H_{m-2}(\Omega (S^{m-1} \times S^{m-1}))$.
  The total homology
$H_*(\Omega (S^{m-1} \times S^{m-1}))$
is isomorphic to the graded commutative
polynomial ring ${\bf Z}[x_1, x_2]$ with
$\deg x_1 = \deg x_2 =m-2$.
Here $x_1$ and $x_2$ correspond to
$\widehat{\alpha}_1$ and $\widehat{\alpha}_2$
respectively.
We have $[x_1, x_2]=0$ where the bracket is
defined by
$[x_1, x_2] = x_1x_2 -(-1)^m x_2x_1$.
This implies the relation
$[X_{ij},  X_{ik} + X_{jk}]=0$.
\end{proof}


The structure of $H_*(\Omega {\rm Conf}_n({\bf R}^m) ; {\bf Z})$
as a Hopf algebra
for $m \geq 3$ was determined by Cohen and
Gitler \cite{CG}.
We denote by
$L_n(m)$ the graded free Lie algebra over $\bf Z$ generated by
$X_{ij}$, $1 \leq i \neq j \leq n$, where the degree for
$X_{ij}$ is $m-2$.
Let $\mathcal T$ denote the ideal of
$L_n(m)$ generated by $X_{ij} - (-1)^m X_{ji}$
together with
\begin{align*}
&[X_{ij},  X_{ik} + X_{jk}],  \quad  i<j<k, \\
&[X_{ij},  X_{kl}], \quad  i, j, k, l \ {\rm distinct}
\end{align*}

We define the Lie algebra ${\mathcal G}_n(m)$ by
$$
{\mathcal G}_n(m) = L_n(m) / {\mathcal T}
$$
and  denote its universal enveloping algebra by
$U{\mathcal G}_n(m)$.
The following theorem is shown inductively by means of
the fibration
$$
\Omega {\rm Conf}_{n+1}({\bf R}^m)
\rightarrow
\Omega {\rm Conf}_n({\bf R}^m).
$$

\begin{thm}{\rm\cite{CG}}\qua
If $m \geq 3$, then we have an isomorphism of graded
Hopf algebras
$$
H_*(\Omega {\rm Conf}_n({\bf R}^m) ; {\bf Z})
\cong
U{\mathcal G}_n(m).
$$
\end{thm}

\section{Chen's iterated integrals}

Our next object is to describe the de Rham cohomology of the
loop space $\Omega {\rm Conf}_n({\bf R}^m)$.
For this purpose we briefly
review Chen's work on the de Rham cohomology
of the loop space.
Let $\omega_1,  \cdots , \omega_q$ be differential forms
on a smooth manifold $M$.
Let $\Delta_q$ be the
$q$-simplex defined by
$$
\Delta_q = \{
(t_1 , \cdots , t_q) \ ; \
0 \leq t_1 \leq \cdots \leq t_q \leq 1 \}.
$$
We have an evaluation map
$$
\phi : \Delta_q \times \Omega M
\rightarrow
M^q
$$
defined by
$$
\phi (t_1 , \cdots t_q ; \gamma) =
(\gamma (t_1) , \cdots \gamma (t_q)).
$$
Chen's iterated integral of the differential forms
$\omega_1 , \cdots , \omega_q$ along the path $\gamma$
is by definition
$$
\int_{\Delta_q} \phi^{*} (
\omega_1 \times \cdots
\times \omega_q)
$$
Following Chen, we denote the above integral
by
$$
\int
\omega_1 \cdots \omega_q.
$$
Let $p$ be the sum of the degrees of $\omega_j$ for
$1\leq j \leq q$.
The iterated integral
$
\int
\omega_1 \cdots \omega_q
$
is considered to be a differential form
of degree $p-q$ on the loop space
$\Omega M$.
Let $B^{p, -q}(M)$ be the vector space over
${\bf R}$ spanned by the iterated integrals of
the form
$
\int
\omega_1 \cdots \omega_q
$
where the sum of the degrees of $\omega_j$,
$1\leq j \leq q$, is equal to $p$.
As a differential form on
the loop space $\Omega M$
$$
d \int
\omega_1 \cdots \omega_q
$$
is expressed up to sign the sum of
$$
\int_{\Delta_q} \phi^{*} d (
\omega_1 \times \cdots
\times \omega_q)
$$
and
$$
\int_{\partial \Delta_q} \phi^{*} (
\omega_1 \times \cdots
\times \omega_q).
$$
Thus we obtain the two differentials
$$
d_1 : B^{p, -q}(M) \rightarrow B^{p+1, -q}(M),
\quad
d_2 : B^{p, -q}(M) \rightarrow B^{p, -q+1}(M).
$$
The direct sum
$$
\bigoplus_{p, q} B^{p, -q}(M)
$$
has a structure of a double complex by the
differentials $d_1$ and $d_2$.
The associated total complex
$B^{\bullet}(M)$ is a subcomplex
of the de Rham complex of the loop space
$\Omega M$.

A basic result due to Chen \cite{Chen} is formulated as follows.
Let us suppose that the manifold $M$ is simply connected.
Then
one has an isomorphism
$$
H^j(\Omega M ; {\bf R}) \cong
H^j(B^{\bullet}(M))
$$
for any $j$.

In the case when $M$ is not simply connected,
the fundamental group of $M$ is related to the
0-dimensional cohomology of the bar complex $B^{\bullet}(M)$
in the following way.
Each element of $H^0(B^{\bullet}(M))$ is represented by
a linear combination of iterated integrals of 1-forms
which is a function on the loop space
$\Omega M$ depending only on the homotopy class of a loop.
Thus we have a natural evaluation map
$$
\pi_1(M, x_0) \times H^0(B^{\bullet}(M))
\rightarrow {\bf R},
$$
which induces a bilinear pairing
$$
{\bf R} \pi_1(M, x_0) \times H^0(B^{\bullet}(M))
\rightarrow {\bf R}.
$$
Here ${\bf R} \pi_1(M, x_0)$ stands for the group ring of
$\pi_1(M, x_0)$ over ${\bf R}$.
Consider the increasing filtration
$$
{\bf R}={\mathcal F}^0 B^{\bullet}(M) \subset
{\mathcal F}^1 B^{\bullet}(M) \subset \cdots \subset
{\mathcal F}^k B^{\bullet}(M) \subset
{\mathcal F}^{k+1} B^{\bullet}(M) \subset \cdots
$$
defined by
$$
{\mathcal F}^k B^{\bullet}(M) =
\bigoplus_{q \leq k} B^{p, -q}(M).
$$
This induces the increasing filtration
$$
{\mathcal F}^k H^0(B^{\bullet}(M))
\subset
{\mathcal F}^{k+1} H^0(B^{\bullet}(M))
$$
on $H^0(B^{\bullet}(M))$.
We denote by ${\mathcal I}$ the augmentation ideal of
${\bf R} \pi_1(M, x_0)$.
It can be shown that for $\omega \in {\mathcal F}^k H^0(B^{\bullet}(M))$
the associated evaluation map
$$
{\bf R} \pi_1(M, x_0)
\rightarrow {\bf R}
$$
factors through ${\mathcal I}^{k+1}$ and we obtain a bilinear map
$$
{\bf R} \pi_1(M, x_0) / {\mathcal I}^{k+1} \times
{\mathcal F}^k B^{\bullet}(M) \rightarrow {\bf R}.
$$
It was shown by Chen (see \cite{Chen}) that
we have an isomorphism
$$
{\mathcal F}^k H^0(B^{\bullet}(M))
\cong
{\rm Hom}_{{\bf R}}
(
{\bf R}\pi_1(M, x_0) / {\mathcal I}^{k+1}, {\bf R}
).
$$

\section{Finite type invariants for braids}

The notion of finite type invariants for braids can be formulated
by means of the group ring in the following way.
We denote by $P_n$ the pure braid group with $n$ strands.
As in the previous section we donote by ${\mathcal I}$ the augmentation
ideal of the group ring
${\bf Z}P_n$.
An invariant $v : P_n \rightarrow {\bf Z}$ is said to be of
order $k$ if the induced map
$v: {\bf Z}P_n \rightarrow {\bf Z}$ factors through
${\mathcal I}^{k+1}$. The set of order $k$ invariants for $P_n$ with
values in
${\bf Z}$ has a structure of a ${\bf Z}$-module and is
identified with
$$
{\rm Hom}_{\bf Z}({\bf Z}P_n/{\mathcal I}^{k+1}, {\bf Z})
$$
and we denote it by
$V_k(P_n)$.
We have a natural inclusion
$V_k(P_n) \subset V_{k+1}(P_n)$ and
we set
$V(P_n) = \bigcup_{k=0}^{\infty}V_k(P_n)$.
We call $V(P_n)$ the space of finite type invariants for
$P_n$ with values in ${\bf Z}$.
The above notion of finite type invariants for $P_n$
is naturally generalized for $B_n$, where $B_n$ stands for
the braid group with $n$ strands.
We consider the ideal ${\mathcal J}$ in the group ring ${\bf
Z}B_n$ generated by
$\sigma_i - \sigma_i^{-1}$,
for a standard system of generators
$\sigma_i$,  $1 \leq i \leq n-1$, for the braid group.
Now an invariant
$v : B_n \rightarrow {\bf Z}$ is defined to be of order $k$
if $v$ vanishies on ${\mathcal J}^{k+1}$.
We denote by $V_k(B_n)$ the set of order $k$ invariants for
$B_n$ with values in ${\bf Z}$
and we set
$V(B_n) = \bigcup_{k=0}^{\infty}V_k(B_n)$.

Considering $P_n$ as the fundamental group of
the configuration space
${\rm Conf}_n({\bf C})$ we readily obtain
the isomorphism
$$
{\mathcal F}^k H^0(B^{\bullet})
\cong
V_k(P_n) \otimes {\bf R}
$$
using Chen's theorem for fundamental groups
(see \cite{Kohno1}).
The above isomorphism provides the expression
of finite type invariants for braids in terms of
iterated integrals of logarithmic forms.
This is a prototype of the Kontsevich integral
in \cite{Kon}.

\begin{figure}[ht!]
\centering
\includegraphics[width=4in]{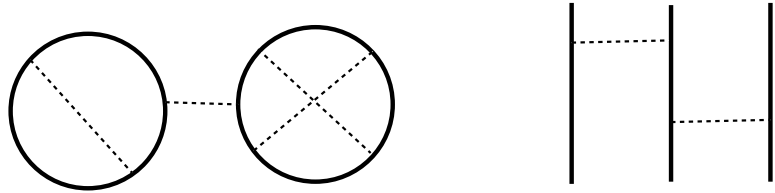}
\caption{A chord diagrams and a horizontal chord diagram}
\label{fig:chord1}
\end{figure}



Let us discuss a relation to the algebra of horizontal chord
diagrams.
First, we recall that
a chord diagram is a collection of
finitely many oriented circles
with finitely many chords attached on them,
regarded up to orientation preserving
diffeomorphisms of the circles (see Figure \ref{fig:chord1}).
Here we assume that the endpoints
of the chords are distinct and lie on
the circles.
Let $I_1 \sqcup \cdots \sqcup I_n$ be the disjoint union of
$n$ unit intervals.
We fix a parameter by a map
$$
p_j : [0, 1] \rightarrow I_j
$$
for each $I_j$, $1\leq j \leq n$.
A horizontal chord diagram on $n$ strands with
$k$ chords is a trivalent graph constructed in the
following way.
We fix $t_1, \cdots, t_k \in [0, 1]$ such that
$0 < t_1 < \cdots < t_k < 1$.
Let
$$
(i_1, j_i), (i_2, j_2), \cdots, (i_k, j_k)
$$
be pairs of distinct integers such that
$1 \leq i_p \leq n$,
$1 \leq j_p \leq n$,
$p=1, 2, \cdots, n$.
We take $k$ copies of unit intervals
$C_1, C_2, \cdots, C_k$ and attach each
$C_{\nu}$ to
$I_1 \sqcup \cdots \sqcup I_n$
in such a way that it starts at
$p_{i_{\nu}}(t_{\nu})$ and ends
at
$p_{j_{\nu}}(t_{\nu})$
for $1 \leq \nu \leq k$.
In this way we obtain a graph
with
$n$ strands $I_1 \sqcup \cdots \sqcup I_n$
and chords $C_1, \cdots, C_k$ attached to them.
Such a graph is called a horizontal
chord diagram on $n$ strands with $k$ chords.
We consider chord diagrams up to orientation
preserving homeomorphism.
Let $D_n^k$ denote the free ${\bf Z}$-module
spanned by  horizontal
chord diagrams on $n$ strands with $k$ chords.
We define ${\Gamma}_{ij}$ as the horizontal chord diagram
on $n$ strands with one chord $C_{ij}$ defined by
the pair $(i, j)$, $1\leq i < j \leq n$.
Then
$$
D_n = \bigoplus_{k \geq 0} D_n^k
$$
has a structure of an algebra
generated by
${\Gamma}_{ij}$,
where the product is defined by the concatenation of
horizontal chord diagrams.

As is the case of the chord diagrams for knots,
we have a natural ${\bf Z}$ module homomorpshim
$$
w : V_k(P_n) \rightarrow {\rm Hom}_{\bf Z}(D_n^k, {\bf Z})
$$
where
$w(v)$ is called the weight system for $v \in V_k(P_n)$.
It can be shown that $w(v)$ vanishes on the ideal
${\mathcal T}$ in $D_n$ generated by
\begin{align*}
&[{\Gamma}_{ij},  {\Gamma}_{ik} + {\Gamma}_{jk}], \quad
[{\Gamma}_{ij} + {\Gamma}_{ik}, {\Gamma}_{jk}] \quad  i<j<k, \\
&[{\Gamma}_{ij},  {\Gamma}_{kl}], \quad  i, j, k, l \ {\rm distinct}.
\end{align*}
We set $A_n = D_n / {\mathcal T}$ and
$A_n^k = D_n^k / D_n^k \cap {\mathcal T}$.
We can show that $w$ induces an isomorphism of ${\bf Z}$ modules
$$
V_k(P_n) / V_{k-1}(P_n) \cong {\rm Hom}_{\bf Z}(A_n^k, {\bf Z})
$$
(see \cite{Kohno1}).

The total homology
$H_*(\Omega {\rm Conf}_n({\bf R}^m); {\bf Z})$ is interpreted as
the algebra of horizontal chord diagrams
on $n$ strands modulo the 4-term relations in the graded sense.
The homology group
$H_{m(k-2)}(\Omega {\rm Conf}_n({\bf R}^m); {\bf Z})$
corresponds to the subspace consisting of
$k$ horizontal chords.
The cohomology
$$
H^{m(k-2)}(\Omega {\rm Conf}_n({\bf R}^m); {\bf Z})
$$
is identified with the space of weight systems on such
horizontal chord diagrams with values in ${\bf Z}$.
We obtain the following theorem (see also \cite{CG}).

\begin{thm}
For $m \geq 3$ we have an isomorphism of ${\bf Z}$-modules
$$
H^{k(m-2)}(\Omega {\rm Conf}_n({\bf R}^m); {\bf Z})
\cong V_k(P_n) / V_{k-1}(P_n).
$$
Let us suppose that $m$ is even.
Then we have an isomorphism of Hopf
algebras
$$
H^*(\Omega {\rm Conf}_n({\bf R}^m); {\bf Z})
\cong V(P_n).
$$
\end{thm}

The above theorem is generalized to
finite type invariants for the full braid group in the
following way.
The symmetric group ${\mathcal S}_n$ acts naturally on
the cohomology ring
$$
H^*(\Omega {\rm Conf}_n({\bf R}^m); {\bf Z})
$$
where the action is induced by the permutations of
the $n$ components of ${\rm Conf}_n({\bf R}^m)$.
Then we have an isomorphism
$$
H^*(\Omega {\rm Conf}_n({\bf R}^m); {\bf Z})
\bullet
{\bf Z}[{\mathcal S}_n]
\cong V(B_n)
$$
where the left hand side stands for the semidirect
product with respect to the above action (see \cite{Kohno2}).
 
\begin{figure}[ht!]
\centering
\includegraphics[width=3.5in]{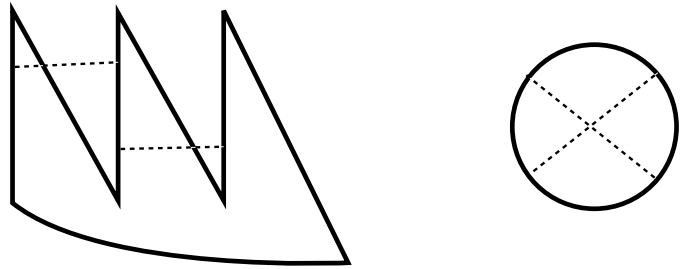}
\caption{Closing a horizontal chord diagram}
\label{fig:chord2}
\end{figure}

Let $\mathcal A$ be the vector space over ${\bf R}$
  spanned
by chord diagrams on a circle
modulo the 4-term relations.
The dual space
$$
{\mathcal A}^* = {\rm Hom}_{\bf R}({\mathcal A},
{\bf R})
$$
is the space of weight systems for Vassiliev invariants
of knots.
We have a surjective map
$$
c : \bigoplus_{n\geq 2}A_n \rightarrow {\mathcal A}
$$
obtained by taking the closure of a horizontal
chord diagram as in Figure \ref{fig:chord2}.
By means of this construction we have the following
theorem.

\begin{thm}
We have an injective homomorphism of
vector spaces
$$
{\mathcal A}^*
\rightarrow
\bigoplus_{n \geq 2} H^*(\Omega {\rm Conf}_n({\bf R}^m);{\bf R}).
$$
for $m\geq 3$.
\end{thm}

Certain cohomology classes in
$\bigoplus_{n \geq 2} H^*(\Omega {\rm Conf}_n({\bf R}^m);{\bf R})$
play a role of
weight systems for Vassiliev invariants.
It would be an interesting problem
to characterize geometrically
such cohomology classes.

\section{Integral representations}

To explain the idea we interpret the Gauss linking number
formula in terms of an integral on the loop space
$\Omega {\rm Conf}_2({\bf R}^3)$.
Let
$$\phi : I \times \Omega {\rm Conf}_2({\bf R}^3)
\rightarrow  {\rm Conf}_2({\bf R}^3)$$
be the evaluation map defined by
$\phi (t, \gamma) = \gamma (t)$.
We pull back $\omega_{12}$ and integrate along the fibre
of the projection map
$p : I \times \Omega {\rm Conf}_2({\bf R}^3)
\rightarrow
\Omega {\rm Conf}_2({\bf R}^3)$
to obtain the 1-form
$$
\int_I \phi^* \omega_{12}
$$
defined on the loop space
$\Omega {\rm Conf}_2({\bf R}^3)$.
This differential form on the loop space
is denoted by
$\int  \omega_{12}$.


A two-component link $L = K_1 \cup K_2$ is
given by a map
$$
f : S^1 \times S^1 \rightarrow {\rm Conf}_2({\bf R}^3).
$$
Let us consider the induced map
$$
\Omega f : \Omega (S^1 \times S^1)
\rightarrow \Omega {\rm Conf}_2({\bf R}^3).
$$
A fundamental homology class of $S^1 \times S^1$
defined by a chain
$\mu : I \times I \rightarrow S^1 \times S^1$
gives in a natural way a 1-chain
$\widehat{\mu} : I \rightarrow \Omega (S^1 \times S^1)$
(see \cite{Chen}).
We set $\widehat{c}_f = \Omega f_* (\widehat{\mu})$
which is a
1-cycle of the loop space
$\Omega {\rm Conf}_2({\bf R}^3)$.
Now we can express the linking number as
$$
lk (K_1, K_2) = \int_{\widehat{c}_f} \int \omega_{12}.
$$


Our next object is to describe the de Rham cohomology of the
loop space
$$
H^*_{DR}(\Omega {\rm Conf}_n({\bf R}^m)).
$$

We see that the
configuration space
${\rm Conf}_n({\bf R}^m)$ is simply connected when
  $m \geq 3$.
Thus we can apply
Chen's theorem to compute the de Rham cohomology of
$\Omega{\rm Conf}_n({\bf R}^m)$.
Let us recall that
the de Rham
cohomology
$H_{DR}^*(\Omega {\rm Conf}_n({\bf R}^m))$
is generated by the
Green forms $\omega_{ij}$, $1\leq i< j \leq n$.
We have the following theorem.

\begin{thm}\label{thm:deRham}
If $m \geq 3$, then
the de Rham cohomology
$$
H_{DR}^{k(m-2)}(\Omega {\rm Conf}_n({\bf R}^m))
$$
is represented by iterated integrals
of the form
$$
\sum a_{i_1j_1 \cdots i_kj_k}
\int
\omega_{i_1j_1} \cdots \omega_{i_kj_k}
+ (\mbox{iterated integrals of length}<k)
$$
where the coefficients $a_{i_1j_1 \cdots i_kj_k} \in {\bf Z}$
satisfy the 4-term relations.
\end{thm}

Here we say that $a_{i_1j_1 \cdots i_kj_k}$ satisfy the
4-term relation if
$X_{i_1j_1 \cdots i_kj_k} \mapsto a_{i_1j_1 \cdots i_kj_k}$
gives a well-defined ${\bf Z}$-module homomorphism
$U{\mathcal G}(m) \rightarrow {\bf Z}$.
In particular, in the case $n=2$,
the de Rham cohomology
$
H_{DR}^*(\Omega {\rm Conf}_2({\bf R}^m))
$
is spanned by the iterated integrals
$$
1, \ \int \omega_{12}, \
\int \omega_{12}\omega_{12},  \ \cdots,
\int \omega_{12} \cdots \omega_{12}, \ \cdots
$$

As an application of the above description of the
de Rham cohomology we give an integral representation
of an order 2 invariant.
Let $L$ be a 3-component link, which is
represented by
$f : S^1 \times S^1 \times S^1 \rightarrow
{\rm Conf}_3({\bf R}^3)$.
As in the case of the Gauss formula for the linking number
we can construct a 2-cycle $\widehat{c}_f$ of
the loop space
$\Omega {\rm Conf}_3({\bf R}^3)$
associated with $f$ in the following way.
We start from a cubic 3-chain
$$
\mu : I^3 \rightarrow S^1\times S^1 \times S^1
$$
corresponding to the fundamental homology class of
the torus $S^1 \times S^1 \times S^1$.
Then $\mu$ gives a 2-chain
$$
\widehat{\mu} : I^2 \rightarrow
\Omega ( S^1 \times S^1 \times S^1).
$$
Composing with
$$
\Omega f :
\Omega ( S^1 \times S^1 \times S^1)
\rightarrow
\Omega {\rm Conf}_3({\bf R}^3)
$$
we obtain a 2-cycle
$\widehat {c}_f$ of the loops space
$\Omega {\rm Conf}_3({\bf R}^3)$.

There is a 3-form $\phi_{123}$ on
${\rm Conf}_3({\bf R}^3)$
such that
the relation
$$
\omega_{12} \wedge \omega_{23} +
\omega_{23} \wedge \omega_{13} +
\omega_{13} \wedge \omega_{12}
= d \phi_{123}
$$
is satisfied.
Now the iterated integral
$$
\int
(\omega_{12} \omega_{23} +
\omega_{23} \omega_{13} +
\omega_{13} \omega_{12})
+ \int \phi_{123}
$$
is a closed 2-form on the loop space
$\Omega {\rm Conf}_3({\bf R}^3)$.
The integral
$$
\int_{\widehat {c}_f}
\left(
\int
(\omega_{12} \omega_{23} +
\omega_{23} \omega_{13} +
\omega_{13} \omega_{12})
+ \int \phi_{123}
\right)
$$
of the above 2-form on the loop space
over the
2-cycle $\widehat{c}_f$
defined above
is an invariant
of $L$. A relation of this invariant to
known invariants will be discussed elsewhere.
We see that
the first term
$$
I_1 = \int_{\widehat {c}_f}
\left(
\int
(\omega_{12} \omega_{23} +
\omega_{23} \omega_{13} +
\omega_{13} \omega_{12})
\right)
$$
is expressed by chord diagrams with 2 chords on
3 circles.
It is important to notice that $I_1$ is not
a link invariant.
A prescription suggested by Chern-Simons perturbation
theory is to compensate this integral by adding
an integral defined by graphs with trivalent vertices
(see \cite{GMM} and \cite{BT}).
 From our point of view, instead of an integral
associated with trivalent graphs, we add
$$
I_2 = \int_{\widehat {c}_f}
\int \phi_{123}
$$
to obtain an integration of a closed 2-form on
the loop space over a 2-cycle.
By suspension we obtain
an $(m-1)$-cycle of
$\Omega {\rm Conf}_3({\bf R}^m)$
associated with
$f : S^1 \times S^1 \times S^1 \rightarrow
{\rm Conf}_3({\bf R}^m)$
for $m>3$ and we have an integral analogous to
the above $I_1 + I_2$ in the case
$m>3$ as well.

Let $\widehat{c_1}$ and $\widehat{c_2}$
be $p$-chain and $q$-chain of the loop space
$\Omega M$.
Then by the composition of loops we can
define a $(p+q)$-chain denoted by
$\widehat{c_1} \cdot \widehat{c_2}$
(see \cite{Chen}).
With this notation the above construction
is generalized in the following way.
Let $L$ be an $n$-component link in ${\bf R}^3$.
We consider an $n$-dimensional torus
$$
T^n = \underbrace{S^1 \times \cdots \times S^1}_{n}.
$$
Then the link $L$ gives a map
$
f : T^n \rightarrow {\rm Conf}_3({\bf R}^3).
$
For a $q$-dimensional subtorus
$T^q \subset T^n$ we have a $q$-cycle
$
\alpha : I^q \rightarrow T^n
$
corresponding to the fundamental homology class of
$T^q$, which gives a $(q-1)$-chain
$$
\widehat{\alpha} : I^{q-1} \rightarrow
\Omega T^n.
$$
Composing with
$$
\Omega f : \Omega T^n \rightarrow
\Omega {\rm Conf}_n({\bf R}^3)
$$
we obtain a $(q-1)$-cycle
$\Omega f_*(\widehat{\alpha})$ of the loop space
$\Omega {\rm Conf}_n({\bf R}^3)$.
Let $\widehat{c}_f$ be a $k$-cycle of
$\Omega {\rm Conf}_n({\bf R}^3)$
represented as the product of cycles of the above type
for any subtorus of $T^n$.
Then we have the following.

\begin{thm}
Let $\omega$ be a closed $k$-form on
$\Omega {\rm Conf}_n({\bf R}^3)$ given by iterated integrals
as in Theorem \ref{thm:deRham}. Then the integral
$$
\int_{\widehat{c}_f} \omega
$$
over a cycle $\widehat{c}_f$ defined
as above
is a link invariant.
\end{thm}

\section{Orbit configuration spaces}

The aim of this section is to give a brief review of
the article \cite{CKX} where we describe the relation
between the homology of the loop spaces of the
orbit configuration associated the action of
Fuchsian groups acting on the upper half plane
and finite type invariants for braids on surfaces.
Let us consider the situation where a group $\Gamma$
acts freely on a space $X$.
We define the orbit configuration space
by
$$
{\rm Conf}_n^{\Gamma}(X)
= \{ (x_1, \cdots, x_n) \in X^n \ ; \
\Gamma x_i \cap \Gamma x_j = \emptyset \
\text{if} \
i \neq j \}
$$
where $\Gamma x$ stands for the orbit of
$x \in X$ with respect to the action of $\Gamma$.

Let $\Gamma$ be a Fuchsian group acting freely on
the upper half plane ${\bf H}$.
The quotient space ${\bf H} / \Gamma$ is an oriented surface
denoted by $\Sigma$.
For $d\geq 1$ we consider the action of
$\Gamma$ on ${\bf H} \times {\bf C}^d$ where
$\Gamma$ acts trivially on ${\bf C}^d$.
We will describe the homology of the loop space of the
orbit configuration space
$$
H_*(\Omega{\rm Conf}_n^{\Gamma}({\bf H} \times {\bf C}^d);
{\bf Z}).
$$
The factor ${\bf C}^d$ appears for the degree shifting.
The reason why we consider the orbit configuration space
rather than the configuration space itself is that
the homology of the loop space of the former one
has a more sensible structure in relation with
finite type invariants for braids on surfaces.

We introduce the notion of horizontal chord diagrams
on $\Sigma$.
First, we recall chord diagrams on surfaces following
\cite{AMR}.
Here we assume that the endpoints
of the chords are distinct and lie on
the circles.
Let $D$ be a chord diagram.
We consider a continuous map
$\gamma : D \rightarrow \Sigma$ and
we denote by $[\gamma]$ its free homotopy class.
We call such pair $(D, [\gamma])$ a chord diagram on $\Sigma$.
We denote by ${\mathcal D}_{\Sigma}$ the vector space over ${\bf R}$
spanned by all chord diagrams on $\Sigma$
and ${\mathcal D}^k_{\Sigma}$ its subspace spanned by
chord diagrams with $k$ chords.
We define ${\mathcal A}(\Sigma)$ to be
the quotient space of ${\mathcal D}_{\Sigma}$
modulo
the 4-term relations.
We refer the reader to \cite{AMR} for a precise definition
of the 4-term relations in this situation.
The chord diagrams on a surface
$\Sigma$ are related to
Vassiliev invariants for links in $\Sigma \times I$ in the
following way.
Let $v$ be an order $k$ invariant for links in
$\Sigma \times I$.
Then the associated weight system $w(v)$ defines a linear form
$w(v) : {\mathcal D}^k_{\Sigma} \rightarrow {\bf R}$
satisfying the 4-term relations.
As in shown in \cite{AMR},
${\mathcal A}(\Sigma)$ has a structure of
a Poisson algebra.

We fix a base point
${\bf x}=(x_1, \cdots, x_n) \in {\rm Conf}_n(\Sigma)$.
The fundamental group
$\pi_1({\rm Conf}_n(\Sigma), {\bf x})$ is the pure braid group
of $\Sigma$ with $n$ strands and is denoted by
$P_n(\Sigma)$.
We have a natural homomorphism
$$
p : P_n(\Sigma) \rightarrow
\bigoplus_{j=1}^n \pi_1(\Sigma, x_j)
$$
and ${\rm Ker} \ p$ is denoted by
$P_n(\Sigma)^0$.
Notice that the direct sum
$\bigoplus_{j=1}^n \pi_1(\Sigma, x_j)$ acts freely on
the orbit configuration space
${\rm Conf}_n^{\Gamma}(H)$ and the quotient space is
the configuration space
${\rm Conf}_n(\Sigma)$.

We denote by $C_n^k$ the set of
horizontal chord diagrams
on $n$ strands with
$k$ chords.
For $\Gamma \in C_n^k$
consider a continuous map
$$
f : \Gamma \rightarrow \Sigma
$$
such that
$$
f(p_i(0))=f(p_i(1))=x_i,  \  1 \leq i \leq n
$$
and denote by
$[f]$ its homotopy class.
Here we consider a homotopy preserving the
base point.
We shall say that such horizontal chord diagram
on $\Sigma$ is
based at
${\bf x}=(x_1, \cdots, x_n)$
We denote by
$D_n^k(\Sigma)$ the free
${\bf Z}$ module
spanned by
pairs
$(\Gamma, [f])$
for $\Gamma \in C_n^k$ and
$
f : \Gamma \rightarrow \Sigma,
$
based at ${\bf x}$.
The subspace of
$D_n^k(\Sigma)$ spanned by
$
f : \Gamma \rightarrow \Sigma,
$
such that each curve
$f(p_i(t)), 0 \leq t \leq 1$, is homotopic to
the point $\{ x_i \}$
is denoted by
$D_n^k(\Sigma)^0$.

We fix a base point $x_0 \in \Sigma$ and
consider the fundamental group
$\pi_1(\Sigma, x_0)$.
Let us fix a path in $\Sigma$
connecting
$x_0$ to $x_j$ and we identify the set of
homotopy classes of paths from
$x_i$ to $x_j$ with
$\pi_1(\Sigma, x_0)$.
For $\gamma \in \pi_1(\Sigma, x_0)$
we consider
$({\Gamma}_{ij}, [f]) \in D_n^k(\Sigma)^0$
such that
$f(C_{ij})$ corresponds to
$\gamma \in \pi_1(\Sigma, x_0)$
by the above identification.
We denote this
$({\Gamma}_{ij}, [f])$
by $X_{ij, \gamma}$.
We see that the direct sum
$$
D_n(\Sigma)^0 = \bigoplus_{k \geq 0}D_n^k(\Sigma)^0
$$
has a structure of an algebra over ${\bf Z}$
where the product is defined by the composition
of chord diagrams.
As an algebra
$D_n^k(\Sigma)^0$
is generated by
$X_{ij, \gamma}$, $1 \leq i \neq j \leq n$,
$\gamma \in \pi_1(\Sigma, x_0)$.
We have the following lemma.

\begin{lem}
The relation $X_{ij, \gamma}=X_{ji, \gamma^{-1}}$
holds for any
$\gamma \in \pi_1(\Sigma, x_0)$.
\end{lem}

The direct sum
$$
D_n(\Sigma) = \bigoplus_{k \geq 0}D_n^k(\Sigma)
$$
has a structure of an algebra as well.
For the subspace
$D_n^0(\Sigma)$ spanned by the chord diagrams with
empty chord, we have a natural injection
$$
\iota_j :
\pi_1(\Sigma, x_j)
\rightarrow
  D_n^0(\Sigma), \ 1 \leq j \leq n
$$
and we have an isomorphism of
${\bf Z}$ algebras
$$
\bigoplus_{j=1}^n
{\bf Z}\pi_1(\Sigma, x_j)
\cong
  D_n^0(\Sigma).
$$
We write $\mu_j$ for $\iota_j(\mu)$
where $\mu$ is an element of  $\pi_1(\Sigma, x_j)$.

The direct sum
$$
\Lambda_n=
\bigoplus_{j=1}^n {\bf Z}
\pi_1(\Sigma, x_j)
$$
acts on
$D_n(\Sigma)^0$ by the conjugation
$$
\Gamma \mapsto \mu_j \Gamma \mu_j^{-1},  \
\Gamma \in D_n(\Sigma)^0,
\mu \in \pi_1(\Sigma, x_j).
$$
We have
the following.

\begin{lem}
With respect to the above action we have
$$
\mu_l X_{ij, \gamma} \mu_l^{-1} = X_{ij, \gamma} \  \mbox{for}
  \ l\neq i, j, \ \
\mu_i X_{ij, \gamma} \mu_i^{-1} = X_{ij, \mu\gamma}.
$$
\end{lem}

Let ${\mathcal I}$ be the ideal of
$D_n(\Sigma)^0$ generated by
\begin{align*}
&[X_{ij, e}, X_{kl, e}], \ \ {i, j, k, l} \ distinct \\
&[X_{ij, e}, X_{jk, e}+X_{ik, e}]  \ \ {i, j, k} \ distinct
\end{align*}
as a $\Lambda_n$ module. We set
$$
{\mathcal A}_n(\Sigma)^0 = D_n(\Sigma)^0 / {\mathcal I},  \ \
{\mathcal A}_n(\Sigma) = D_n(\Sigma) / {\mathcal I}.
$$
We have an action of $\Lambda_n$ on
${\mathcal A}_n(\Sigma)^0$ by conjugation and the semidirect
product
${\mathcal A}_n(\Sigma)^0 \bullet \Lambda_n$ with respect to
this action is isomorphic to
${\mathcal A}_n(\Sigma)$.

\begin{lem}\label{lem:4t}
We have
the following relations in ${\mathcal A}_n(\Sigma)^0$.
\begin{align*}
&[X_{ij, \gamma}, X_{kl, \delta}]=0, \ \ {i, j, k, l} \ distinct \\
&[X_{ij, \gamma}, X_{jk, \delta}+X_{ik, \gamma\delta}]=0  \ \ {i, j, k} \ 
distinct
\end{align*}
\end{lem}

The above 4-term relations appear naturally when we
consider finite type invariants for
$P_n(\Sigma)$.
We denote by
${V}_k(P_n(\Sigma))$ the
${\bf Z}$ module of
order $k$  invariants of
$P_n(\Sigma)$ with values in ${\bf Z}$.
For $v \in {V}_k(P_n(\Sigma))$ the associated weight system
$w(v)$ defines a ${\bf Z}$-module homomorphism
$$
w(v) :
D^k_n(\Sigma) \rightarrow
{\bf Z}
$$
satisfying the 4-term relations in Lemma \ref{lem:4t}.
We have an injective homomorphism of
${\bf Z}$-modules
$$
{V}_k(P_n(\Sigma))/ {V}_{k-1}(P_n(\Sigma))
\rightarrow
{\rm Hom}({\mathcal A}_n^k(\Sigma), {\bf Z})
$$
where ${\mathcal A}_n^k(\Sigma)$ denotes the submodule of
${\mathcal A}_n(\Sigma)$ spanned by chord diagrams with
$k$ chords.

The algebra  ${\mathcal A}_n(\Sigma)^0$ and ${\mathcal A}_n(\Sigma)$
were introduced in \cite{gp} by the above generators and
relations.
In this article, they constructed an injective homomorphism
$$
{\bf Z}P_n(\Sigma) \rightarrow
\widehat{{\mathcal A}_n(\Sigma)} =
\prod_{k \geq 0}{\mathcal A}_n^k(\Sigma)
$$
as ${\bf Z}$ modules.

The homology of the loop space of our orbit configuration
space is related to the algebra of horizontal chord diagrams
on $\Sigma = {\bf H} / \Gamma$
  in the following way.

\begin{thm}{\rm\cite{CKX}}\qua
We have an isomorphism of Hopf algebras
$$
H_*(\Omega{\rm Conf}_n^{\Gamma}({\bf H} \times {\bf C}^d);
{\bf Z})
\cong
{\mathcal A}_n(\Sigma)^0
$$
for any $d \geq 1$.
\end{thm}

\Addresses

\end{document}